\documentclass{article}

\usepackage{amsmath,amssymb,amsthm}
\usepackage{mathrsfs}
\usepackage{amscd}
\usepackage{enumerate}
\theoremstyle{plain}
\newtheorem{theorem}{Theorem}[section]
\theoremstyle{remark}
\newtheorem{remark}[theorem]{Remark}

\newtheorem{definition}[theorem]{Definition}
\theoremstyle{plain}
\newtheorem{corollary}[theorem]{Corollary}

\numberwithin{equation}{section}

\def\N{{\mathbb N}}
\def\Z{{\mathbb Z}}

\def\R{{\mathbb R}}

\usepackage{color}

\definecolor{gr}{rgb}   {0.,   0.8,   0. } 
\definecolor{bl}{rgb}   {0.,   0.5,   1. } 
\definecolor{mg}{rgb}   {0.7,  0.,    0.7} 

\newcommand{\Bk}{\color{black}}

\title{The maximal regularity operator on tent spaces}

\author{Pascal Auscher, Sylvie Monniaux, Pierre Portal}

\date{En l'honneur des 60 ans de Michel Pierre}

\begin{document}

\maketitle

\begin{abstract}
Recently, Auscher and Axelsson gave a new approach to non-smooth boundary value 
problems with $L^{2}$ data, that relies on some appropriate weighted maximal regularity 
estimates. As part of the development of the corresponding $L^{p}$ theory, we prove here 
the relevant weighted maximal estimates in tent spaces $T^{p,2}$ for $p$ in a certain open 
range. We also study the case $p=\infty$.
\end{abstract}

\section{Introduction}
Let $-L$ be a densely defined closed linear operator acting on $L^{2}(\R^{n})$ and generating 
a bounded analytic semigroup $(e^{-tL})_{t \geq 0}$.
We consider the maximal regularity operator defined by
$$
\mathcal{M}_{L}f(t,x) = \int \limits _{0} ^{t} Le^{-(t-s)L}f(s,.)(x)ds,
$$
for functions $f \in C_{c}(\R_{+} \times \R^{n})$. 
The boundedness of this operator on $L^{2}(\R_{+}\times \R^{n})$ was established by 
de\,Simon in \cite{desimon}. The $L^{p}(\R_{+}\times \R^{n})$ case, for $1<p<\infty$, turned 
out, however, to be much more difficult. In \cite{kl}, Kalton and Lancien proved that 
$\mathcal{M}_{L}$ could fail to be bounded on $L^p$ as soon as $p \neq 2$. The necessary 
and sufficient assumption for $L^p$ boundedness was then found by Weis \cite{w} to be a 
vector-valued strengthening of analyticity, called R-analyticity. 
As many differential operators $L$ turn out to generate R-analytic semigroups, the $L^p$ 
boundedness of $\mathcal{M}_{L}$ has subsequently been successfully used in a variety of 
PDE situations (see \cite{kuw} for a survey). 

Recently, maximal regularity was used in a different manner as an important tool in \cite{aa}, where a new 
approach to boundary 
value problems with $L^{2}$ data for divergence form elliptic systems on Lipschitz domains, 
is developed. More precisely, in \cite{aa}, the authors establish and use the boundedness of 
$\mathcal{M}_{L}$ on weighted spaces $L^{2}(\R_{+}\times \R^{n}; t^{\beta}dt dx)$, for certain 
values of $\beta \in \R$, under the additional assumption that  $L$ has bounded holomorphic functional calculus on $L^2(\R^n)$. This additional assumption  was removed in \cite[Theorem~1.3]{aa2}. Here is the version when specializing the Hilbert space to be $L^2(\R^n)$.
\begin{theorem} \label{AA} With $L$ as above,    $\mathcal{M}_{L}$ extends to a bounded operator on $ L^{2}(\R_{+}\times \R^{n}; t^{\beta}dt dx)$ for all $\beta \in (-\infty, 1)$.
\end{theorem}
 The use of these weighted spaces is common in the study of boundary value 
problems, where they are seen as variants of the tent space $T^{2,2}$ which occurs for $\beta=-1$,  introduced by Coifman, 
Meyer and Stein in \cite{cms}. For $p \neq 2$, the corresponding spaces are weighted versions 
of the tent spaces $T^{p,2}$, which are defined, for parameters $\beta \in \R$ and 
$m \in \N$, as the completion of $C_{c}(\R_{+} \times \R^{n})$ with respect to 
$$
\|g\|_{T^{p,2,m}(t^{\beta}dt dy)} = \left( \int \limits _{\R^{n}} 
\Bigl( \int \limits _{0} ^{\infty} \int \limits _{\R^{n}} \frac{1_{B(x,t^{\frac{1}{m}})}(y)}{t^{\frac{n}{m}}} \,
\bigl|g(t,y)\bigr|^{2} t^{\beta}dydt \Bigr)^{\frac{p}{2}}dx \right)^{\frac{1}{p}},
$$
the classical case corresponding to $\beta = -1$, $m=1$, and being denoted simply by $T^{p,2}$.
The parameter $m$ is used to allow various homogeneities, and thus to make these spaces 
relevant in the study of differential operators $L$ of order $m$.
To develop an analogue of \cite{aa} for $L^p$ data, we need, among many other estimates yet to be proved, boundedness results for the 
maximal operator $\mathcal{M}_{L}$ on these tent spaces. This is the purpose of this note. 
Another motivation  is 
well-posedness of non-autonomous Cauchy problems for operators with varying domains, which will be presented elsewhere. In the latter case,  $\mathcal{M}_{L}$ can be seen as a model of the evolution operators involved.  However, as $\mathcal{M}_{L}$ is an important operator on its own, we thought interesting to present this special case alone. 

In Section \ref{sec:result} we state and prove the adequate boundedness results. The proof  is based on recent results and methods developed in \cite{hnp}, building on ideas 
from \cite{amr} and \cite{hmp}.
In Section \ref{sec:tools} we recall the relevant material from \cite{hnp}.

\section{Tools}
\label{sec:tools}

When dealing with tent spaces, the key estimate needed is a change of aperture formula, i.e., 
a comparison between the $T^{p,2}$ norm and the norm
$$
\|g\|_{T^{p,2}_\alpha} := \left( \int \limits _{\R^{n}} \Bigl( \int \limits _{0} ^{\infty} 
\int \limits _{\R^{n}} \frac{1_{B(x,\alpha t)}(y)}{t^{n}}\, \bigl|g(t,y)\bigr|^{2} 
\frac{dy dt}{t}\Bigr)^{\frac{p}{2}}dx \right)^{\frac{1}{p}},
$$
for some parameter $\alpha>0$.
Such a result was first established in \cite{cms}, building on similar estimates in \cite{fs}, 
and analogues have since been developed in various contexts. Here we use the following 
version given in \cite[Theorem 4.3]{hnp}.

\begin{theorem}\label{thm:angle}
Let $1<p<\infty$ and $\alpha \geq 1$. There exists a constant $C>0$ such that, for all 
$f\in T^{p,2}$, 
\begin{equation*}
   \| f\|_{T^{p,2}}\leq \| f\|_{T_{\alpha}^{p,2}}\leq C(1+\log\alpha)\alpha^{n/\tau}
   \| f\|_{T^{p,2}},
\end{equation*}
where $\tau = \min(p,2)$ and $C$ depends only on $n$ and $p$. \footnote{only on pour \'eviter les confusions}
\end{theorem}

Theorem \ref{thm:angle} is actually a special case of the Banach space valued result obtained 
in \cite{hnp}.
Note, however, that it improves the power of $\alpha$ appearing in the inequality from the 
$n$ given in \cite{cms} to $\frac{n}{\tau}$.
This is crucial in what follows, and has been shown to be optimal in \cite{hnp}. 

Applying this to $(t,y) \mapsto t^{\frac{m(\beta+1)}{2}}f(t^{m},y)$ instead of $f$, we also have 
the weighted result, where
$$
\|g\|_{T^{p,2,m} _{\alpha}(t^{\beta}dt dy)} = \left( \int \limits _{\R^{n}} \Bigl( \int \limits _{0} ^{\infty} 
\int \limits _{\R^{n}} \frac{1_{B(x,\alpha t^{\frac{1}{m}})}(y)}{t^{\frac{n}{m}}} \,
\bigl|g(t,y)\bigr|^{2} t^{\beta}dydt \Bigl)^{\frac{p}{2}}dx \right)^{\frac{1}{p}}.
$$

\begin{corollary}\label{cor:angle}
Let $1<p<\infty$, $m \in \N$, $\alpha \geq 1$, and $\beta \in \R$. There exists a constant 
$C>0$ such that, for all $f\in T^{p,2,m}(t^{\beta}dt dy)$, 
\begin{equation*}
   \| f\|_{T^{p,2,m}(t^{\beta}dt dy)}\leq \| f\|_{T_{\alpha}^{p,2,m}(t^{\beta}dt dy)}
   \leq C(1+\log\alpha)\alpha^{n/\tau}
   \| f\|_{T^{p,2,m}(t^{\beta}dt dy)},
\end{equation*}
where $\tau = \min(p,2)$ and $C$ depends only on $n$ and $p$.
\end{corollary}

To take advantage of this result, one needs to deal with families of operators, that behave 
nicely with respect to tent norms.
As pointed out in \cite{hnp}, this does not mean considering R-bounded families (which means 
R-analytic semigroups when one considers $(tLe^{-tL})_{t\geq0}$) as in the 
$L^{p}(\R_{+}\times \R^{n})$ case, but tent bounded ones, i.e. families of operators with the 
following $L^{2}$ off-diagonal decay, also known as Gaffney-Davies estimates.

\begin{definition}
A family of bounded linear operators $(T_{t})_{t \geq 0} \subset B(L^{2}(\R^{n}))$ is said to 
satisfy off-diagonal estimates of order $M$, with homogeneity $m$, if, for all Borel sets 
$E,F \subset \R^{n}$, all $t>0$, and all $f \in L^{2}(\R^{n})$:
$$
\|1_{E}T_{t}1_{F}f\|_{2} \lesssim \Big(1+\frac{dist(E,F)^{m}}{t}\Big)^{-M}\|1_{F}f\|_{2}.
$$
In what follows $\|\cdot \|_{2}$ denotes the norm in $L^2(\R^n)$.
\end{definition}

As proven, for instance, in \cite{ahlmt}, many differential operators of order $m$, such as 
(for $m=2$) divergence form elliptic operators with bounded measurable complex coefficients, 
are such that $(tLe^{-tL})_{t\geq0}$ satisfies off-diagonal estimates of any order, with 
homogeneity $m$. This condition can, in fact, be seen as a replacement for the classical 
gaussian kernel estimates satisfied in the case of more regular coefficients.

\section{Results}
\label{sec:result}

\begin{theorem}
\label{thm:main}
Let $m \in \N$, $\beta \in (-\infty,1)$, 
 $p \in \bigl(\frac{2n}{n+m(1-\beta)},\infty\bigr) \cap (1,\infty)$, and $\tau=\min(p,2)$. If 
 $(tLe^{-tL})_{t \geq 0}$ satisfies off-diagonal estimates of order 
$M>\frac{n}{m\tau}$, with homogeneity $m$,  then $\mathcal{M}_{L}$ extends to a bounded operator on 
$T^{p,2,m}(t^{\beta}dt dy)$.
\end{theorem}

\begin{proof}
The proof is very much inspired by similar estimates in \cite{amr} and \cite{hnp}.
Let $f \in {\mathscr{C}}_{c}(\R_{+}\times \R^{n})$. 
Given $(t,x) \in \R_{+} \times \R^{n}$, and $j \in \Z_{+}$, we consider 
$$
C_{j}(x,t) = 
\begin{cases} B(x,t) \; \text{if} \; j=0, \\
B(x,2^{j}t)\backslash B(x,2^{j-1}t) \; \text{otherwise.} \end{cases}
$$
We write $\|\mathcal{M}_{L}f\|_{T^{p,2}} \leq \sum \limits _{k=1} ^{\infty} \sum \limits 
_{j=0} ^{\infty} I_{k,j} + \sum \limits _{j=0} ^{\infty} J_{j}$ where
\begin{eqnarray*}
I_{k,j} &=& \left( \int \limits _{\R^{n}} \Bigl( \int \limits _{0} ^{\infty} \int \limits _{\R^{n}} 
\frac{1_{B(x,t^{\frac{1}{m}})}(y)}{t^{\frac{n}{m}}} \,
\Bigl|\int \limits _{2^{-k-1}t} ^{2^{-k}t} Le^{-(t-s)L}(1_{C_{j}(x,4t^{\frac{1}{m}})}f(s,.))(y)ds\Bigr|^{2} \, 
{t^\beta dy\,dt}\Bigr)^{\frac{p}{2}}dx \right)^{\frac{1}{p}}, 
\\[4pt]
J_{j} &=& \left( \int \limits _{\R^{n}} \Bigr( \int \limits _{0} ^{\infty} \int \limits _{\R^{n}} 
\frac{1_{B(x,t^{\frac{1}{m}})}(y)}{t^{\frac{n}{m}}} \,
\Bigl|\int \limits _{\frac{t}{2}} ^{t} Le^{-(t-s)L}(1_{C_{j}(x,4s^{\frac{1}{m}})}f(s,.))(y)ds\Bigr|^{2} \, 
{t^\beta dy\,dt}\Bigr)^{\frac{p}{2}}dx \right)^{\frac{1}{p}}.
\end{eqnarray*}
Fixing $j \geq 0$, $k \geq 1$ we first estimate $I_{k,j}$ as follows. For  fixed $x\in {\mathbb{R}}^{n}$,
\begin{eqnarray*}
\label{est:low}
&&\int_{0}^{\infty} \int_{B(x,t^{\frac{1}{m}})} \Bigl| \int_{2^{-k-1}t}^{2^{-k}t} Le^{-(t-s)L}
(1_{C_{j}(x,4t^{\frac{1}{m}})}f(s,\cdot))(y)\,ds\Bigl|^{2} t^{\beta-\frac{n}{m}} dy\,dt 
\\[4pt]
&\leq&
\int_{0} ^{\infty} \int_{B(x,t^{\frac{1}{m}})} \Bigl(\int_{2^{-k-1}t}^{2^{-k}t} 
\Bigl|(t-s)Le^{-(t-s)L}(1_{C_{j}(x,4t^{\frac{1}{m}})}f(s,\cdot))(y)\Bigr|\,\frac{ds}{t-s}\Bigr)^{2} 
t^{\beta-\frac{n}{m}} dy\,dt 
\\[4pt]
&\lesssim&
\int_{0}^{\infty} \int_{2^{-k-1}t}^{2^{-k}t}  2^{-k}t \Bigl(\int_{B(x,t^{\frac{1}{m}})}  
\bigl|(t-s)Le^{-(t-s)L}(1_{C_{j}(x,4t^{\frac{1}{m}})}f(s,\cdot))(y)\bigr|^{2} dy \Bigr)\,
t^{\beta-\frac{n}{m}-2} ds\,dt 
\\[4pt]
&\lesssim&
\int_{0}^{\infty} \int_{2^{-k-1}t}^{2^{-k}t} 
2^{-k}\Big(1+\frac{2^{jm}t}{t-s}\Big)^{-2M}\,
\bigl\|1_{B(x,2^{j+2}t^{\frac{1}{m}})}f(s,\cdot)\bigr\|_{2}^{2} \,t^{\beta - \frac{n}{m} -1}ds\,dt
\\[4pt]
&\lesssim&
2^{-k}2^{-2jmM}\int_{0}^{\infty} 
\Bigl(\int_{2^ks}^{2^{k+1}s}t^{\beta-\frac{n}{m}-1}dt\Bigr)
\bigl\|1_{B(x,2^{j+\frac{k}{m}+3}s^{\frac{1}{m}})}f(s,\cdot)\bigr\|_{2}^{2} \,ds
\\[4pt]
&\lesssim&
2^{-k(\frac{n}{m}+1-\beta)}2^{-2jmM}\int_{0}^{\infty} 
\bigl\|1_{B(x,2^{j+\frac{k}{m}+3}s^{\frac{1}{m}})}f(s,\cdot)\bigr\|_{2}^{2} 
\,s^{\beta-\frac{n}{m}}ds.
\end{eqnarray*}
In the second inequality, we use Cauchy-Schwarz inequality for the integral with respect to $t$, 
the fact that $t-s\sim t$ for $s\in\cup_{k\ge1}[2^{-k-1}t,2^{-k}t] \subset[0,\frac{t}{2}]$ and Fubini's 
theorem to exchange the integral in $t$ and the integral in $y$. The next inequality follows from 
the off-diagonal estimate verified by $(t-s)Le^{-(t-s)L}$ and again the fact that $t-s\sim t$. 
By Corollary \ref{cor:angle} this gives
$$
I_{k,j} \lesssim (j+k)2^{-k(\frac1 2({\frac{n}{m}+1-\beta})-\frac{n}{m\tau})}2^{-j(mM-\frac{n}{\tau})}
\|f\|_{T^{p,2,m}(t^{\beta}dt dy)},
$$
where $\tau=\min(p,2)$.
It follows \Bk that 
$\sum\limits_{k=1}^{\infty} \sum\limits_{j=0}^{\infty} I_{k,j} \lesssim 
\|f\|_{T^{p,2,m}(t^{\beta}dt dy)}$
since $M>\frac{n}{m\tau}$ and ${\frac{n}{m}+1-\beta} > \frac{2n}{m\tau}$ (Note that for 
$p\ge 2$, this requires $\beta<1$). 

We now turn to $J_{0}$ and remark that
$J_{0} \leq \bigl( \int  _{\R^{n}} J_{0}(x)^{\frac{p}{2}}dx \bigr)^{\frac{1}{p}}$, where
\begin{equation*}
J_{0}(x)=\int \limits _{0} ^{\infty} \int \limits _{\R^{n}} 
\Bigl|\int \limits _{\frac t 2} ^{{t}{}} Le^{-(t-s)L}(g(s,\cdot)(y)ds\Bigr|^{2}\, 
t^{\beta-\frac{n}{m}}dy\,dt
\end{equation*}
with $g(s,y)=1_{B(x,4s^{\frac{1}{m}})}(y)f(s,y)$.
The inside integral can be rewritten as 
$$
\mathcal{M}_{L}g(t,\cdot) -e^{-\frac{t}{2} L} \mathcal{M}_{L}g(\frac t 2,\cdot).
$$ 
As $\mathcal{M}_{L}$ is bounded on $L^2(\R_{+}\times \R^{n}; t^{\beta-\frac{n}{m}}dy dt)$ by 
Theorem~\ref{AA} and $(e^{-tL})_{t\ge0}$ is uniformly bounded on $L^{2}(\mathbb{R}^n)$, 
we get 
$$
J_0(x) \lesssim \int_{0}^{\infty} \bigl\|1_{B(x,4s^{\frac{1}{m}})}f(s,\cdot)\bigr\|_{2}^{2}\,
s^{\beta-\frac{n}{m}}ds.
$$

We finally turn to $J_{j}$, for $j \ge 1$. For fixed  $x\in \R^{n}$, 
\begin{eqnarray*}
&&\int \limits _{0} ^{\infty}  \int\limits_{\R^{n}} 1_{B(x,t^{\frac{1}{m}})}(y)  
\Bigl| \int\limits_{\frac{t}{2}}^{t} Le^{-(t-s)L}(1_{C_{j}(x,4s^{\frac{1}{m}})}f(s,.))(y)ds\Bigr|^{2} 
t^{\beta-\frac{n}{m}}dy\,dt 
\\[4pt]
&\leq& 
\int \limits _{0} ^{\infty} \int \limits _{\R^{n}} 1_{B(x,t^{\frac{1}{m}})}(y) 
\Bigl(\int \limits _{\frac{t}{2}} ^{t} \bigl|(t-s)Le^{-(t-s)L}(1_{C_{j}(x,4s^{\frac{1}{m}})}f(s,.))(y)\bigr|
\frac{ds}{t-s}\Bigr)^{2} t^{\beta-\frac{n}{m}}dy\,dt 
\\[4pt]
&\lesssim& 
\int \limits_{0} ^{\infty} \int \limits _{\R^{n}} 1_{B(x,t^{\frac{1}{m}})}(y) 
\int \limits_{\frac{t}{2}} ^{t} 
\bigl|(t-s)Le^{-(t-s)L}(1_{C_{j}(x,4s^{\frac{1}{m}})}f(s,.))(y)\bigr|^{2}\frac{ds}{(t-s)^{2}} 
t^{\beta-\frac{n}{m}+1}dy\,dt 
\\[4pt]
&\lesssim& 
\int \limits _{0} ^{\infty} \int \limits _{\frac{t}{2}} ^{t} (t-s)^{-2} 
\Bigl(1+\frac{2^{jm}t}{t-s}\Bigr)^{-2M}
\bigl\|1_{B(x,2^{j+2}s^{\frac{1}{m}})}f(s,.)\bigr\|_{2} ^{2}\, s^{\beta-\frac{n}{m}+1}ds\,dt
\\[4pt]
&\lesssim& 
2^{-jm(2M-2)}\int \limits _{0} ^{\infty} \left(\int \limits _{s} ^{2s} s(t-s)^{-2}
\Bigl(1+\frac{2^{jm}t}{t-s}\Bigr)^{-2}dt\right)
\bigl\|1_{B(x,2^{j+2}s^{\frac{1}{m}})}f(s,.)\bigr\|_{2} ^{2}\,  s^{\beta-\frac{n}{m}}ds
\\[4pt]
&\lesssim& 
2^{-2jmM}\int \limits _{0} ^{\infty} 
\bigl\|1_{B(x,2^{j+2}s^{\frac{1}{m}})}f(s,.)\bigr\|_{2} ^{2}\,
s^{\beta-\frac{n}{m}}ds,
\end{eqnarray*}
where we have used Cauchy-Schwarz inequality in the second inequality, the off-diagonal 
estimates and the fact that $s \leq t$ in the third, Fubini's theorem and the fact that 
$s \geq \frac{t}{2}$ in the fourth, and the change of variable $\sigma= \frac{t}{t-s}$ in the last.
An application of Corollary \ref{cor:angle}, then gives
$$ 
J_{j} \lesssim 2^{-jmM}j2^{j\frac{n}{\tau}}\|f\|_{T^{p,2,m}(t^{\beta}dtdy)} = 
j 2^{-j(mM-\frac{n}{\tau})}\|f\|_{T^{p,2,m}(t^{\beta}dtdy)}, 
$$
and the proof is concluded by summing the estimates.
\end{proof}

An end-point result holds for $p=\infty$. In this context the appropriate tent space consists of 
functions such that $|g(t,x)|^{2}\frac{dx dt}{t}$ is a Carleson measure, and is defined as the 
completion of the space ${\mathscr{C}}_{c}(\R_{+}\times \R^{n})$ with respect to
$$
\|g\|_{T^{\infty,2}}^2=  
\underset{(x,r) \in \R^{n} \times \R_{+}}{\sup} 
r^{-n} \int \limits _{B(x,r)} \int \limits _{0} ^{r} |g(t,x)|^{2} \frac{dx dt}{t}.
$$
We also consider the weighted version defined by 
$$
\|g\|_{T^{\infty,2,m}(t^{\beta}dtdy)} ^{2} := 
\underset{(x,r) \in \R^{n} \times \R_{+}}{\sup} 
r^{-\frac{n}{m}} \int \limits _{B(x,r^{\frac{1}{m}})} \int \limits _{0} ^{r} |g(t,x)|^{2} t^{\beta}dx dt.
$$

\begin{theorem}
\label{thm:infini}
Let $m \in \N$, and $\beta \in (-\infty,1)$.
If $(tLe^{-tL})_{t \geq 0}$ satisfies off-diagonal estimates of order $M>\frac{n}{2m}$, with 
homogeneity $m$, then $\mathcal{M}_{L}$ 
extends to a bounded operator on $T^{\infty,2,m}(t^{\beta}dtdy)$.
\end{theorem}

\begin{proof} 
Pick a ball $B(z,r^{\frac{1}{m}})$. Let 
$$
I^2=\int \limits _{B(z,r^{\frac{1}{m}})} 
\int \limits _{0} ^{r} |(\mathcal{M}_{L}f)(t,x)|^{2} t^{\beta}dx dt.
$$
We want to show that $I^{2} \lesssim r^{\frac{n}{m}} \|f\|_{T^{\infty,2}(t^{\beta}dtdy)}^2$. 
We set 
$$
I_{j}^2=
\int \limits _{B(x,r^{\frac{1}{m}})} \int \limits _{0} ^{r} |(\mathcal{M}_{L}f_{j})(t,x)|^{2} t^{\beta}dx dt
$$ 
where $f_{j}(s,x)=f(s,x) 1_{C_{j}(z,4r^{\frac{1}{m}})}(x)1_{(0,r)}(s)$ for $j\ge 0$. Thus by 
Minkowsky inequality, $I \le \sum I_{j}$. For $I_{0}$ we use again Theorem~\ref{AA} 
which implies that $ \mathcal{M}_{L}$ is bounded on $L^2(\R_{+}\times \R^{n}, t^{\beta}dxdt)$. 
Thus 
$$
I_{0}^2\lesssim  \int \limits _{B(z,4r^{\frac{1}{m}})} \int \limits _{0} ^{r} |f(t,x)|^{2} t^{\beta}dx dt 
\lesssim r^{\frac{n}{m}} \|f\|_{T^{\infty,2,m}(t^{\beta}dtdy)}^2.
$$
Next, for $j \neq 0$, we proceed as in the proof of Theorem \ref{thm:main} to obtain
\begin{eqnarray*}
I_{j}^{2} 
&\lesssim& \sum \limits _{k=1} ^{\infty} \int \limits _{0} ^{r} \int \limits _{2^{-k-1}t} ^{2^{-k}t}
2^{-k}t \Bigl(1+\frac{2^{jm}r}{t-s}\Bigr)^{-2M} \|f_{j}(s,.)\|_{L^{2}} ^{2} t^{\beta-2}ds\,dt
\\[4pt]
&&+ \int \limits _{0} ^{r} \int \limits _{\frac{t}{2}} ^{t} t(t-s)^{-2} \Bigl(1+\frac{2^{jm}r}{t-s}\Bigr)^{-2M}
\|f_{j}(s,.)\|_{L^{2}} ^{2} t^{\beta}ds\,dt.
\end{eqnarray*}
Exchanging the order of integration, and using the fact that $t \sim t-s$ in the first part and 
that $t \sim s$ in the second, we have the following.
\begin{eqnarray*}
I_{j}^{2} &\lesssim& 
\sum \limits _{k=1} ^{\infty} 2^{-k}2^{-2jmM}r^{-2M}
\int \limits _{0} ^{2^{-k}r} \int \limits _{2^{k}s} ^{2^{k+1}s} t^{\beta+2M-1}
\|f_{j}(s,.)\|_{L^{2}} ^{2} dtds 
\\[4pt]
&& + \int \limits _{0} ^{r} \int \limits _{s} ^{2s} r(t-s)^{-2}
\Bigl(1+\frac{2^{jm}r}{t-s}\Bigr)^{-2M}\|f_{j}(s,.)\|_{L^{2}} ^{2} s^{\beta}dtds 
\\[4pt]
&\lesssim&
\sum \limits _{k=1} ^{\infty} 2^{-k}2^{-2jmM}
\int \limits _{0} ^{2^{-k}r} (2^{k}s)^{\beta} \|f_{j}(s,.)\|_{L^{2}} ^{2} ds 
+ \int \limits _{0} ^{r} \int \limits _{1} ^{\infty} \bigl(1+2^{jm}\sigma\bigr)^{-2M}
\|f_{j}(s,.)\|_{L^{2}} ^{2} s^{\beta}d\sigma ds
\\[4pt]
&\lesssim& 2^{-2jmM} \int \limits _{0} ^{r}  \|f_{j}(s,.)\|_{L^{2}} ^{2} s^{\beta}ds,
\end{eqnarray*}
where we used $\beta<1$. We thus have 
$$
I_{j}^{2} \lesssim 2^{-2jmM}(2^{j}r^{\frac{1}{m}})^{n}\|f\|^{2}_{T^{\infty,2,m}(t^{\beta} dtdy)},
$$
and the condition $M>\frac{n}{2m}$ allows us to sum these estimates. 
\end{proof}

\begin{remark}
Assuming off-diagonal estimates, instead of kernel estimates, allows to deal with differential 
operators $L$ with rough coefficients.
The harmonic analytic objects associated with $L$ then fall outside the Calder\'on-Zygmund 
class, and it is common (see for instance \cite{memoir}) for their boundedness range to be a 
proper subset of $(1,\infty)$. Here, our range $(\frac{2n}{n+m(1-\beta)},\infty]$ includes 
$[2,\infty]$ as $\beta<1$, which is consistent with \cite{aa}. In the case of classical tent spaces, 
i.e., $m=1$ and $\beta=-1$, it is the range $(2_{*},\infty]$, where $2_{*}$ denotes the Sobolev 
exponent $\frac{2n}{n+2}$. We do not know, however, if this range is optimal.
\end{remark}

\begin{remark}
Theorem \ref{thm:infini} is a maximal regularity result for parabolic Carleson measure norms.
This is quite natural from the point of view of non-linear parabolic PDE (where maximal regularity 
is often used), and such norm have, actually, already been used in the context of Navier-Stokes
equations in \cite{kt}, and, subsequently, for some geometric non-linear PDE in \cite{KoLa}.
Theorem \ref{thm:main} is also reminiscent of Krylov's Littlewood-Paley estimates \cite{Kry}, and of their 
recent far-reaching generalization in \cite{nvw}.
In fact, the methods and results from \cite{hnp}, on which this paper relies, use the same circle of 
ideas (R-boundedness, Kalton-Weis $\gamma$ multiplier theorem...) as \cite{nvw}. The combination 
of these ideas into a ``conical square function" approach to stochastic maximal regularity will be 
the subject of a forthcoming paper.
\end{remark}

{\flushleft{\sc Pascal Auscher}}\\
Univ. Paris-Sud, laboratoire de Math\'ematiques, UMR 8628, F-91405 {\sc Orsay}; CNRS, F-91405 {\sc Orsay}. \\
{\tt pascal.auscher@math.u-psud.fr}\\

{\flushleft{\sc Sylvie Monniaux}}\\
LATP-UMR 6632, FST Saint-J\'er\^ome - Case Cour A, 
Univ. Paul C\'ezanne, F-13397 {\sc Marseille} C\'edex 20.\\
{\tt sylvie.monniaux@univ-cezanne.fr}\\

{\flushleft{\sc Pierre Portal}}\\
Permanent Address:\\
Universit\'e Lille 1, Laboratoire Paul Painlev\'e, F-59655 {\sc Villeneuve d'Ascq}.\\
Current Address:\\
Australian National University, Mathematical Sciences Institute, John Dedman Building, 
Acton ACT 0200, Australia.\\
{\tt pierre.portal@math.univ-lille1.fr}\\

\end{document}